\documentclass{amsart}
\usepackage{amssymb}
\usepackage{epsfig}
\usepackage{amsfonts}
\usepackage{epsfig}

\newcommand{\Hom}{{\mathrm{Hom}}}

\theoremstyle{plain}
\newtheorem{theorem}{Theorem}
\newtheorem{lemma}{Lemma}

\newtheorem{corollary}[theorem]{Corollary}

\theoremstyle{definition}

\theoremstyle{remark}
\newtheorem{remark}[theorem]{Remark}

\def\C{\mathbb C}
\def\O{\mathcal O}

\begin{document}

\title[An upper bound for the lower central series
quotients]{An upper bound for the lower central series
quotients of a free associative algebra}

\author{Galyna Dobrovolska}
\address{Department of Mathematics, Massachusetts Institute of
Technology, Cambridge, MA 02139, U.S.A.}
\email{galyna@mit.edu}

\author{Pavel Etingof}
\address{Department of Mathematics, Massachusetts Institute of
Technology, Cambridge, MA 02139, U.S.A.}
\email{etingof@math.mit.edu}

\maketitle

\begin{abstract}

Feigin and Shoikhet conjectured in \cite{FS} that successive quotients $B_m(A_n)$ of the lower central series filtration of a free associative algebra $A_n$ have polynomial growth. In this paper we give a proof of this conjecture, using the structure of a representation of $W_n$, the Lie algebra of polynomial vector fields on $\C^n$, on $B_m(A_n)$ which was defined in \cite{FS}. Moreover, we show that the number of squares in a Young diagram $D$ corresponding to an irreducible $W_n$-module in the Jordan-H\"older series of $B_m(A_n)$ is bounded above by the integer $(m-1)^2+2[\frac{n-2}{2}](m-1)$, which allows us to confirm the structure of $B_3(A_3)$ conjectured in \cite{FS}.

\end{abstract}

\section{introduction}

Let $A:=A_n$ be the free associative algebra  with $n$ generators
over $\mathbb{C}$. Consider its lower central series as a Lie algebra,
i.e., the Lie ideals $L_m(A)\subset A$ defined recursively by 
$L_1(A)=A, L_i(A)=[A,L_{i-1}(A)]$, and the corresponding
associated graded Lie algebra $B(A):=\oplus_{i\ge 1}B_i(A)$, 
where $B_i(A)=L_i(A)/L_{i+1}(A)$. 

Giving each of the generators of $A$ degree $1$, 
we define a grading on $A$ and hence on each of the spaces
$B_i(A)$. It is an interesting (and, in general, 
unsolved) problem to determine the Hilbert series of $B_i(A)$ for each $i$ 
with respect to this grading. For $i=1$, this is easy since 
$B_1(A)$ is the space of cyclic words in $n$ letters.
It particular, one can easily see that $B_1(A)$ has exponential growth.  
So at first sight one might expect that the spaces $B_i(A)$ 
for $i>1$ also have exponential growth. However, computer
experiments performed by Eric Rains in 2005 suggested 
that, to the contrary, the spaces $B_i(A)$ for $i>1$ should have polynomial
growth. 

This phenomenon was studied systematically by Feigin and Shoikhet
in the paper \cite{FS}. The first important observation of Feigin
and Shoikhet is that the image $Z$ of $A[[A,A],A]A$ in $B_1(A)=A/[A,A]$ is central in
the Lie algebra $B(A)$, and therefore, it is natural to define the
space $\bar B_1(A)=B_1(A)/Z$ and the graded Lie algebra $\bar B(A)=\bar
B_1(A)\oplus \oplus_{i\ge 2}B_i(A)$. Their second important
observation is that the graded Lie algebra $\bar B(A)$ carries a
natural grading-preserving action of the Lie algebra 
$W_n$ of polynomial vector fields in $n$ variables
(while there is no such action on $B_1(A)$), and the $W_n$-modules
$\bar B_1(A)$ and $B_i(A)$, $i\ge 2$, admit a Jordan-H\"older
series whose simple composition factors are irreducible $W_n$-modules
${\mathcal F}_D$ of tensor fields corresponding to Young diagrams $D$ with at most
$n$ rows.

Feigin and Shoikhet in \cite{FS} determined the exact structure of $\bar B_1(A)$
and $B_2(A)$ as representations of $W_n$, showing that they have
polynomial growth, and conjectured that
for any $m$, $B_m(A)$ has polynomial growth, i.e.,
the dimension $B_m(A_n)[\ell]$, the degree $\ell$ part of
$B_m(A_n)$, grows like $c_{mn}\ell^{n-1}$, where $c_{mn}$ is a constant.

In this paper we give a proof of this conjecture. In fact, we
show that the $W_n$-modules $B_m(A_n)$ have finite length,
which implies the conjecture. More specifically, we give an explicit
upper bound on the number of squares in a Young diagram $D$ if 
the corresponding tensor field module ${\mathcal F}_D$ over $W_n$
occurs as a composition factor in $B_m(A_n)$. This bound not only implies that $B_m(A_n)$
has finite length, but also, when combined with the computation by Eric Rains
(see \cite{FS}), allows us to confirm the conjectural
structure of $B_3(A_3)$ (as conjectured in \cite{FS}, it turns
out to be isomorphic to ${\mathcal F}_{2,1,0}$).

\begin{remark} In this paper we work over the ground field $\mathbb{C}$
following \cite{FS}, but the discussion carries over without changes to
an arbitrary field of characteristic zero.
\end{remark}

{\bf Acknowledgments.} P.E. is grateful to B. Shoikhet for useful discussions.
The work of P.E. was  partially supported by the NSF grant
 DMS-0504847.

\section{Results}

\subsection{Tensor field modules over $W_n$}
The simple $W_n$-modules ${\mathcal F}_D$ are defined in the following way. Let $D$ be a Young diagram with at most $n$ rows. Then ${\mathcal F}_D$ is the irreducible $W_n$-submodule of the space
of polynomial tensor fields on $\C^n$ of type $D$. In other words, we take a Young diagram $D$ and the corresponding irreducible $\mathfrak{gl}(n,\mathbb{C})$-module $F_D$. Then we extend the action of $\mathfrak{gl}(n,\mathbb{C})$ on $F_D$ to the action of $W_n^0$, the subalgebra of $W_n$ of vector fields vanishing at the origin, so that quadratic and higher vector fields act by zero. Then, if $D$ has more than one column, we let ${\mathcal{F}}_D=\Hom_{U(W_n^0)}(U(W_n),F_D)$ be the $W_{n}$-module coinduced from a $\mathfrak{gl}(n,\mathbb{C})$-module $F_D$, which is known to be irreducible. If $D$ has one column, we take ${\mathcal{F}}_D$ to be the irreducible $W_n$-submodule of the coinduced module $\Hom_{U(W_n^0)}(U(W_n),F_D)$ (which is known to be unique in this case).

It is known that if $D$ has more than one column, then ${\mathcal F}_D$ is the whole space of tensor fields, otherwise (if $D$ has one column of $k$ squares) ${\mathcal F}_D$ is the space of
closed polynomial $k$-forms on $\mathbb C^n$. If
$(k_1,\ldots,k_n)$ is a partition (possibly ending with some
zeroes) corresponding to a Young diagram $D$, we will write
${\mathcal F}_{(k_1,\ldots,k_n)}$ to denote the module ${\mathcal
F}_D$. For references about 
modules ${\mathcal{F}}_D$, see \cite{FF} or \cite{F}; for reference about Schur modules $F_D$, see \cite{Ful}.

\subsection{The main result}
The main result of this paper is the following theorem, proved in
the next section.

\begin{theorem}\label{main} For $m\geq 3,\ n\geq 2$
and $\mathcal{F}_D$ in the Jordan-H\"older  series of $B_m(A_n)$,
we have the following estimate on the size (i.e., the number of
squares) of the Young diagram $D$:
$$
|D|\leq (m-1)^2+2[\frac{n-2}{2}](m-1)
$$
(where $[x]$ denotes the integer part of $x$).
\end{theorem}

\begin{corollary} For all $m\ge 2,n\ge 2$, $B_m(A_n)$ has finite length
as a $W_n$-module. In particular, $\dim B_m(A_n)[\ell]\sim
c_{mn}\ell^{n-1}$, as $\ell\to \infty$, where $c_{mn}>0$ is a constant.
\end{corollary}

\begin{proof} It suffices to consider the case $m\ge 3$, as
the case $m=2$ is considered in \cite{FS}.
The first statement follows from Theorem \ref{main} and the fact that
$\dim B_m(A_n)[\ell]<\infty$. The second statement follows from
the first statement.
\end{proof}

\begin{corollary}
$B_3(A_3)={\mathcal F}_{2,1,0}$.
\end{corollary}

\begin{proof}
According to the MAGMA computation by Eric Rains (see \cite{FS},
formula (17) for $H_3(u,t)$), the characters of both sides agree
up to degree 6 inclusively (see \cite{FS}, section 4.2). So if in
addition to ${\mathcal F}_{2,1,0}$, $B_3(A_3)$ had included another
${\mathcal F}_D$, then $D$ would have to have at least 7 squares.
But according to Theorem \ref{main}, $D$ can have at most 4
squares. Thus, there is no additional constituents, and we are
done.
\end{proof}

\begin{remark}
In a similar way, one can rederive the results of \cite{DKM}
on the structure of $B_3(A_2)$ and $B_4(A_2)$. Namely,
for $B_4(A_2)$, the bound for the size of $D$ given by
Theorem \ref{main} is 9, while the Hilbert series of $B_4(A_2)$
is computed by Eric Rains using MAGMA up to degree 9 (see
\cite{FS}).
\end{remark}

\section{Proof of Theorem \ref{main}}

\subsection{The map $\psi$ and its kernel}

Consider the polynomial rings \linebreak $\O_{mn}:=\C[x_{ij}]_{1\leq i \leq
n, 1\leq j \leq m}$. Below we write $A:=A_n$. Let
$\Omega^{even}(\Bbb C^n)$ be the space of even polynomial
differential forms, and set $\Omega^{mn}:=
\Omega^{even}(\C^n)^{\otimes m}$. As an $\O_{mn}$-module,
$\Omega^{mn} \cong \O_{mn} \otimes \Lambda^{mn}$, where
$\Lambda^{mn}:=\Lambda^{even}(\C^n)^{\otimes m}$.

Recall that by a result of \cite{FS}, we have a surjective map
$\xi:\Omega^{even}(\C^n)\to \bar B_1(A)$, which descends to an isomorphism
$\Omega^{even}(\C^n)/\Omega^{even}_{exact}(\C^n) \cong
\bar{B}_1(A)$, which we will also denote by $\xi$.
For any $m\ge 2$, this isomorphism gives rise to a map $\psi:
\Omega^{mn} \rightarrow B_m(A)$ defined as a composition of two maps
$$
\Omega^{mn}= \Omega^{even}(\C^n)^{\otimes m} \rightarrow
\bar{B}_1(A)^{\otimes m} \rightarrow B_m(A),
$$
$\omega_1 \otimes \dots \otimes \omega_m \mapsto b_1
\otimes \dots \otimes b_m \mapsto [[b_1,b_2],\ldots, b_m]$,
where $b_i:=\xi(\omega_i)$. The map $\psi$ is surjective, since
by the results of \cite{FS}, the direct sum $\bar B_1(A)\oplus
\bigoplus_{m\ge 2}B_m(A)$ is a graded Lie algebra generated in
degree 1.

Recall also from \cite{FS} that we have a natural isomorphism
$\eta: \Omega^{even >0}_{exact}(\Bbb C^n)\to B_2(A)$,
and upon identification by $\xi,\eta$, the bracket map
$\bar B_1(A)^{\otimes 2}\to B_2(A)$ reduces to the map
$\omega_1\otimes\omega_2\to d\omega_1\wedge d\omega_2$.

Therefore, we can factor the map $\psi$ as
$$
\Omega^{even}(\C^n)^{\otimes 2} \otimes
\Omega^{even}(\C^n)^{\otimes (m-2)} \rightarrow B_2(A) \otimes
\bar{B}_1(A)^{\otimes (m-2)} \rightarrow B_m(A)
$$
using the map
$\eta \circ \phi: \Omega^{even}(\C^n)^{\otimes 2} \rightarrow B_2(A),$ where
$\phi: \omega_1 \otimes \omega_2 \mapsto d \omega_1 \wedge d\omega_2.$

 \begin{lemma} \label{computation} For any $i,j \in [1,n]$ and
 forms $\omega_1,\omega_2 \in \Omega^{even}(\C^n)$, the element
\linebreak $(x_{i1}-x_{i2})(x_{j1}-x_{j2}) \omega_1 \otimes \omega_2$ is in
 the kernel of the map $\phi$. Consequently, the submodule
 $(x_{i1}-x_{i2})(x_{j1}-x_{j2}) \Omega^{mn}$ is in the kernel of
 $\psi$.
\end{lemma}

\begin{proof}
We have
$$
\phi((x_{i1}-x_{i2}) \omega_1 \otimes \omega_2)
= \{ d(x_{i1} \omega_1)\wedge d \omega_2 - d \omega_1
\wedge d(x_{i2} \omega_2) \} |_{x_{i1}=x_{i2}=x_i}
$$
$$
= dx_i \wedge \omega_1 \wedge d \omega_2+ x_i d \omega_1 \wedge d
\omega_2 - d \omega_1 \wedge dx_i \wedge \omega_2 - d \omega_1
\wedge x_i d \omega_2
$$
(since $\omega_1$ is even)
$$
=dx_i \wedge (\omega_1 \wedge d \omega_2 + d \omega_1 \wedge
\omega_2).
$$
$$
= dx_i \wedge d(\omega_1 \wedge \omega_2).
$$
Using this, we compute
$$
\phi((x_{i1}-x_{i2})(x_{j1}-x_{j2}) \omega_1 \otimes \omega_2) =
$$
$$
= \phi((x_{j1}-x_{j2})(x_{i1}\omega_1 \otimes \omega_2) -
(x_{j1}-x_{j2})(\omega_1 \otimes x_{i2} \omega_2))
$$
$$
=dx_j \wedge d(x_i \omega_1 \wedge \omega_2) - dx_j \wedge d(
\omega_1 \wedge x_i \omega_2) = 0.
$$
Therefore, $\phi$ maps $(x_{i1}-x_{i2})(x_{j1}-x_{j2})
\Omega^{even} (\C^n)^{\otimes 2}$ to zero. In particular,
this implies that $(x_{i1}-x_{i2})(x_{j1}-x_{j2}) \Omega^{mn}$
belongs to the kernel of the map $\psi$.
\end{proof}

\begin{lemma}\label{comid} Let $L$ be a Lie algebra, and $b_i\in
L$, $i=1,...,m$.
For any $k=1,...,m$, the bracket $[[b_1,b_2],\ldots,
b_m]$ is a linear combination of the brackets of the form
$[[b_k,b_{l_1}], \dots, b_{l_{m-1}} ]$ where
$(l_1,\ldots,l_{m-1})$ is a permutation of $[1,m]\setminus\{k\}.$
\end{lemma}

\begin{proof}
This lemma is well known, but we will give a proof, as it is very
short. The proof is by induction in $k$. Indeed, for
$k=1,2$ the statement is true. To go from $k-1$ to $k$, we notice that by
the Jacobi identity
$$
[[b_1,b_2],\ldots, b_{k-1}], b_k] = [[b_1,b_2],\ldots, b_k],
b_{k-1}] + [[b_1,b_2],\ldots], [b_{k-1}, b_k]].
$$
Putting $b_{k-1}^{'}:=b_k$ we have that the first bracket on the
RHS can be expressed as a linear combination of brackets of the
form $[[b_k,b_{l_1}], \dots, b_{l_{k-1}} ]$ by the induction
assumption. Similarly, in the second bracket on the RHS we put
$b_{k-1}^{''}:=[b_{k-1},b_k]$ and use
the induction assumption to express the second bracket as a
combination of brackets of the form $[[b_k,b_{k-1}],\ldots].$
Bracketing the LHS and RHS with $b_{k+1},\ldots, b_m$, we obtain the desired result.
\end{proof}

Let $I$ be the ideal in $\O_{mn}$ generated by the polynomials
$$
\prod_{s=1}^{k-1} (x_{i_s s}- x_{i_s k})(x_{j_s s}- x_{j_s k}).
$$
where $2\leq k \leq m$ and arbitrary $i_s,j_s \in [1,n]$ for $s\in [1,k]$.

\begin{lemma} \label{I(Omega) vanishes} The submodule
$I \Omega^{mn}$ of $\Omega^{mn}$ is $W_n$-invariant, and
is contained in the kernel of the map $\psi$.
\end{lemma}

\begin{proof}

First we show that $I \Omega^{mn}$ is stable under the action of
of $W_n$ on $\Omega^{mn}$. For this, it suffices to show that
the ideal $I$ is $W_n$-invariant. Let $f \frac{\partial}{\partial y_m}$ be
some vector field in $W_n$. Then we have:
$$
f \frac{\partial}{\partial y_m} x_{ij}= \left\{ \begin{array}{ll}
0 & \textrm{if $k \neq j$}\\
f_j & \textrm{if $m=i$}
\end{array} \right.
$$
where $f_j$ denotes $f$ where instead of the variables $y_l$ we
substitute $x_{lj}$. Therefore,
$$
(f \frac{\partial}{\partial y_m})
\prod_{s=1}^{k-1} (x_{i_s s}- x_{i_s k})(x_{j_s s}- x_{j_s k})$$
$$= \prod_{s=1}^{k-1} (x_{i_s s}- x_{i_s k})(x_{j_s s}- x_{j_s k})
\sum_{s=1}^{k-1}(\frac{f_s - f_k}{x_{i_s s} - x_{i_s k}}
+\frac{f_s - f_k}{x_{j_s s} - x_{j_s k}})
$$
The last factor is clearly a polynomial, so the right hand side
belongs to $I$. This implies that $I$ is $W_n$-invariant, as desired.

Next we prove that the element $\prod_{s=1}^{k-1} (x_{i_s
s}-x_{i_s k})(x_{j_s s}-x_{j_s k})\omega_1 \otimes \dots \otimes
\omega_m$ goes to zero under the map $\psi.$ Let $b_i=\xi(\omega_i)$. By Lemma \ref{comid}, in
$B_m(A)$ we have
$$
[[b_1,b_2],\ldots, b_m] = \sum_{\sigma}
c_\sigma[[b_k,b_{\sigma(1)}],\ldots, b_{\sigma(m-1)}],
$$
where $\sigma$ is a bijection from $[1,m-1]$ to $[1,m]\setminus
\lbrace{k\rbrace}$.
Therefore, we have that
$$
\omega_1 \otimes \dots \otimes \omega_m - \sum_{\sigma}
c_\sigma(\omega_k\otimes \omega_{\sigma(1)})\otimes \dots \otimes \omega_{\sigma(m-1)}
$$
is in the kernel of $\psi.$

But by Lemma \ref{computation} $(x_{i_{\sigma(1)}
\sigma(1)}-x_{i_{\sigma(1)} k})(x_{j_{\sigma(1)}
\sigma(1)}-x_{j_{\sigma(1)} k})
(\omega_k\otimes\omega_{\sigma(1)})
\otimes \dots \otimes \omega_{\sigma(m-1)}$ is in the kernel
of $\psi.$

Therefore, $I \Omega^{mn}$ is in the kernel of $\psi.$
\end{proof}

Set $r:=2[\frac{n-2}{2}](m-1)$. Note that the space $\Lambda(\Bbb
C^n)$ is equipped with the natural grading by rank of exterior
forms; hence so is the space $\Lambda^{mn}$.
Consider the space
$\Lambda^{mn}_{>r}$ spanned by homogeneous elements of
degree $>r$ in $\Lambda^{mn}$.

\begin{lemma} \label{big forms vanish}
The submodule $\O_{mn} \otimes \Lambda^{mn}_{>r}$ of $\Omega^{mn}$
 is $W_n$-invariant, and is contained in the kernel of $\psi$.
\end{lemma}

\begin{proof}
Since vector fields act as Lie derivative on the components
of $\Omega^{mn}$, they leave rank unchanged, so
$\O_{mn} \otimes \Lambda_{> r}$ is stable under the $W_n$-action on
$\Omega^{mn}$.

Now we show that this submodule is annihilated by $\psi$. Elements of
$\psi(\O_{mn} \otimes \Lambda^{mn}_{>r})$ are
linear combinations of elements of the form $[[b_1, b_2] \ldots b_m]$,
$b_i=\xi(\omega_i)$,
where $\sum \mathrm{rk} \omega_i>2[\frac{n-2}{2}](m-1)$. By
Lemma \ref{comid}, for every $k\in [1,m]$ we have
$[[b_1, b_2] \ldots b_m]= \sum_{\sigma} c_\sigma[[b_k,
b_{\sigma(1)}] \ldots b_{\sigma(m-1)}]$.

Pick $k$ such that $\omega_k$ has maximal rank among
$\omega_1,\ldots,\omega_m$. Then $\mathrm{rk}
\omega_k+\mathrm{rk} \omega_{\sigma(1)}\geq \frac{1}{m-1}\sum
\mathrm{rk}\omega_i >2[\frac{n-2}{2}]$. But
$\mathrm{rk}\omega_k+\mathrm{rk}\omega_{i_{1}}$ is even so
$\mathrm{rk}\omega_k+\mathrm{rk}\omega_{i_{1}}\geq
2([\frac{n-2}{2}]+1)$. So we have
$$
\mathrm{rk}(d\omega_k \wedge d\omega_{\sigma(1)})=
\mathrm{rk}\omega_k+\mathrm{rk}\omega_{\sigma(1)}+2 \geq
2([\frac{n-2}{2}]+2)=2[\frac{n+2}{2}]>n.
$$
Thus $d\omega_k \wedge d\omega_{\sigma(1)}=0$ in $\Omega^{even}(\C^n).$

But $[b_k,b_{\sigma(1)}]=\eta(d\omega_k \wedge d\omega_{\sigma(1)})$ in $B_2(A)$.
So $[b_k,b_{\sigma(1)}]=0$, and thus we have
$[[b_1, b_2] \ldots b_m]=0$.
Therefore the submodule $\O_{mn} \otimes \Lambda^{mn}_{>r}$ is in the kernel of $\psi$.
\end{proof}

\subsection{The structure of $O_{mn}/I$ as a graded space.}

Note that the algebra $\O_{mn}$ has a natural grading, which
assigns degree 1 to each generator, and $I$ is a graded ideal.
Thus, $\O_{mn}/I$ is a graded algebra.

\begin{lemma} \label{O/I} We have an isomorphism of graded vector
spaces
$$
\O_{mn}/I \cong \C[y_1,\ldots,y_n] \otimes
\bigotimes_{k=2}^{m} \mathrm{Sym}^{2k-3} (\C+\sum_{i=1}^n \C
y_i),
$$
where $y_i$ have degree $1$.
\end{lemma}

\begin{proof} Denote by $\mathrm{in}(I)$ the set of
initial (highest) terms of polynomials in $I$ with respect to the
lexicographic monomial order $x_{ij}>x_{kl}$ iff $j>l$ or $j=l
\text{ and } i>k$. By a theorem about Gr\"obner bases (Theorem
15.3 in \cite{E}) monomials not in $\mathrm{in}(I)$ form a vector
space basis of $\O_{mn}/I$. By the form of the generators of $I$
we see that $\mathrm{in}(I)$ is generated by the monomials
$\prod_{s=1}^{2k-2}x_{i_s k}$, where $k\geq 2,\ i_s \in [1,n]$ is
arbitrary. So monomials not in $\mathrm{in}(I)$ are precisely the
monomials in $\C[x_{11},\ldots,x_{n1}] \otimes
\bigotimes_{k=2}^{m} \mathrm{Sym}^{2k-3} (\C+\sum_{i=1}^n \C
x_{ik}).$
\end{proof}

\subsection{Proof of Theorem \ref{main}}

For any graded space $M$ define its character
(the Hilbert series of $M$) $\mathrm{char} M= \sum_{a} \dim
M[a]t^a$, where $M[a]$ is the $a^{th}$
graded piece of $M$. For instance,
$\mathrm{char}{\mathcal F}_D=\frac{P_D(t)}{(1-t)^n}$, where
$P_D(t)$ is a polynomial of degree $n$ if $D$ has one column,
and $P_D(t)=N_Dt^{|D|}$ if $D$ has more than one column,
where $N_D$ is the dimension of the irreducible representation of
$GL_n$ corresponding to $D$.

The space $\Omega^{mn}=\O_{mn}\otimes \Lambda^{mn}$
has a natural grading coming from the grading on the
factors, and we have a surjective map of graded vector spaces
$\psi:\Omega^{mn}
\twoheadrightarrow B_m(A_n)$. By Lemma \ref{I(Omega) vanishes},
$I\Omega^{mn}=I\O_{mn} \otimes \Lambda^{mn}$ is in the kernel of
$\psi$, and by Lemma \ref{big forms vanish}, $\O_{mn} \otimes
\Lambda^{mn}_{>r}$ is in the kernel of $\psi$. So we have a
surjective map
$$(\O_{mn} \otimes \Lambda^{mn})/(I\O_{mn} \otimes
\Lambda^{mn}+\O_{mn} \otimes \Lambda^{mn}_{>r})
\twoheadrightarrow B_m(A_n).$$

Thus we have a surjective map
$$
\hat\psi: (\O_{mn}/I)\otimes (\Lambda^{mn}/ \Lambda^{mn}_{>r})
\twoheadrightarrow B_m(A_n),
$$
which is a homomorphism of $W_n$-modules. We also notice that the space $\Lambda^{mn}/ \Lambda^{mn}_{>r}$ is naturally identified with
$\Lambda^{mn}_{\leq r}$.

Let $\mu(D)$ be the multiplicity of ${\mathcal F}_D$ in
$(\O_{mn}/I)\otimes (\Lambda^{mn}/ \Lambda^{mn}_{>r})$, and $\nu(D)$ be the multiplicity of
${\mathcal F}_D$ in $B_m(A_n)$. Since $\hat\psi$ is surjective,
$\nu(D)\le \mu(D)$.
Also, we have
\begin{equation}\label{equu}
\sum_D \mu(D){\rm char}{\mathcal F}_D=
{\rm char}(\O_{mn}/I){\rm char}(\Lambda^{mn}_{\le r}).
\end{equation}

 From Lemma
\ref{O/I} we have
$$
\mathrm{char} \O_{mn}/I = \mathrm{char} \C[y_1,\ldots,y_n]
\prod_{k=2}^{m} \mathrm{char}\ \mathrm{Sym}^{2k-3}
(\C+\sum_{i=1}^n \C y_i)
$$
$$
=\frac{1}{(1-t)^n} \prod_{k=2}^{m} (\sum_{a=0}^{2k-3}
{{a+n-1} \choose n}t^a)
$$
$$
=\frac{1}{(1-t)^n} (c_0 t^{(m-1)^2}+\mathrm{l.o.t.}),
$$
where $c_0\neq 0$ is a constant and l.o.t. denote monomials of smaller degree than $(m-1)^2$.
The character of $\Lambda^{mn}_{\leq r}$ is a polynomial of degree $r$, $C_0 t^r +\mathrm{l.o.t.}$

Thus, multiplying (\ref{equu}) by $(1-t)^n$, we obtain
$$
\sum_D\mu(D)P_D(t)=Q(t),
$$
where $Q(t)$ is a polynomial of degree $(m-1)^2+r$.

Now, for $m\geq 3,n\ge 2$
$$
n\le 2(n-1)\leq 4[n/2]=4+4[(n-2)/2]\le (m-1)^2+r.
$$
So, we get
$$
\sum_D \mu(D)N_D t^{|D|}=Q_*(t),
$$
where $Q_*$ is another polynomial of degree $(m-1)^2+r$.
Hence $\mu(D)=0$ for $|D|>(m-1)^2+r$ and hence
$\nu(D)=0$ for $|D|>(m-1)^2+r$, as desired.

\end{document}